\documentclass{article}
\usepackage{amsmath,amssymb,amsfonts}

\begin{document}
	\title{Example of  cube slices that  are  not zonoids }
	\author{R.Anantharaman}
	\date{March 15, 2018}
\maketitle

\textbf{To the memories of Som   Naimpally and Joe Diestel}

\begin{abstract}
	
	Let Q be the unit cube in  $ R^ n$  centered    at the Origin O and H a hyperplane   through   O.The intersection is     called a central Cube slice  and  its study was    initiated by Hadwiger,  Henesley  and Vaaler , continued by Ball   and others. A zonoid is the range of a non atomic vector   measure into $R^n$ . In this paper, when $n=4$ we give  examples of  non -zonoid cube slices.   Let H: x + y +z +t =0 ; the slice has triangle faces and is not a zonoid. This contrasts with   a result in$ R^3$,where it  follows from a classical Theorem  due to Herz   and  Lindenstrauss  that  every  central  cube slice is a zonoid( zonotope).  We also  give   nontrivial examples in which the   slice is a  zonoid.  For  ex. let  H : ax + y + z+t=0 with   a$ > $1.  If  a${\ge}$ 3 , the slice is a zonotope. Otherwise it has faces that are trapeziums    or  pentagons and is not a zonoid.We also give other examples of the like nature.
	
\end{abstract}                                                  

\textbf{  Published :Math Japon.  SCMJ e-2018  whole number 31, pp   315-326}

\textbf{New Reference  Ethan  Bolker, at end of paper  JULY 2023}

\section{Introduction}\label{S:intro}

\textbf{1.1  Slices Zonoids ,Zonotopes}

Let us recall the result from [3]:-- Let $Q^n$ =Q =unit cube in $ R ^n$ centered at Origin O;
ie. Q = $ ( (x_ k) :  | x_k|  \le \dfrac{1}{2}  $   .

Let H be a     vector subspace of dimension  n-1, ie. a plane thru the Origin  with equation : H= (  x=  ( $x _k $ ) with  x  .a    =0  ) for a( non zero){\bf} vector  a   in $ R ^n$. The intersection of   H andQ will be called central slice or, “ slice”.  Following [ 3 ] we denote by$ | A|$  the  appropriate volume /area   of the measurable set A  ,and assume n ${\ge }2$. As other examples let us note  the papers  [7 ], [ 8],[ 13]initial to this subject  , and the surveys [5], [10] [14]  ]that treats  many related topicsWenote the p th powerof $ L^ p $ norm of the “sinc “ function in [3] : for (p$\ge2$ ):

\begin{equation}\label{E:cong1}
  I_p =\dfrac{1}{\pi} \int_{R}\dfrac{|sin t|^p}{|t||^p}dt ,
\end{equation}

An upperbound for this is:-
$\dfrac{\sqrt2}{\sqrt p}\ $, with equality iff $p=2$. The lower bound is  assumed  by H:  $x_k$= 0   and upper   only if n=2   and H: $  x + y$=0  or with $ x - y$=0

Let us mention that Valler[13] considers  concepts   of analytic interest;  his results not only prove lower bound but also apply to Minkowski's Theorem on Linear foirms.
  
We note that there are   alsoresults on sections  by   centralplanes of dimension k( see [14] TH 1.2, 1.3 p 154),  alsodue toBall. We  treat only the case  k= n-1.

This  estimate  is in [3 ]; see also [10, Ch1].    The  proof of this estimatein [3, p468] is with  " direct  " and  uses  only  elementary methods . The one in [10] uses  Fourier methods.This integral $I_p$  has found use in wavelets[ 11 ]

For our needs we use the more precise values also from ( [3] Lemma3) below  ,see  eq(9),  (10).
In [3] this is  derived ,first using Characteristic functions( = Fourier Transform)   then the “ standard  Inverse  Fourier Formula ”  .

As pointed out by an anonymous   referee (of another paper)-- see Acknowledgements --this  $I_4$  is in  the classic, [ 12]    ( also in [9]); see [10] for  many    related deeper results  
However  we   use the   formula    from  [3] for vol of   slice of  cube . Our  interest  is more in the sliceitself  .With n=4   in Sec 3we give example ofa a ( central)slice  that has a triangle face, and  isnot nota zonoid ( ${" face"}$ ) defined below).On the other hand,in Sec 4    we  give examples of slices   that are zonoids, and  othersthat have a pentagon or trapezium face and so are not.

\textbf{Notation  and preliminaries :} We write an element of $R^4$ as  (x, y, z, t) and use a, b, c, d as coefficients. Below  we  avoid the   case when  H is paralll to acoordinate  hyperplane;inthis case the slice is a Cube of lower dimension and so a zonoid

Let a hyperplane be H : ax  + by + cz + dt =0 .As  in [ 8 ], we may assume that no coefficient is zero, and next they are all positive. Further we may assumethat $a{\ge} b{\ge}c{\ge}d$
and  then by dividing by d,  that H: ax + by + cz + t =0  with $ a{\ge}b{\ge}c{\ge}1  $.   In all examples of non zonoids we consider the   equation [ t=  -1/2] to get a  Face  that is atriangle,trapezium or pentagon( disqualifying slice fombeig a zonoid : see  beginning of Sec 3).

In ex 3.1  we consider the case  when a=b=c(= 1); and as mentioned above show that the slice has triangular faces and so not a zonoid(" face" defined below).The sections  of this sliceby planes[ t= -c] with $0<c<1/2$  are hexagons .These tend to the triangle face as c  tends to  1/2 .
 We may  feel that "  cube slices in $R^4$  are never  nontrivial zonoids".  Hencein ex 4.1 we consider  H: a x  + y + z + t =0   with $ a >$ 1 Now the slice is a zonoid  if $  a {\ge}$3  and is a paralletope;in the contrary cases the slice has  pentagon faces and isnot a zonoid . In Ex4.2 we consider H:a( x +y)+ z +t=0 ; the slice is not a zonoid onaccount of trapezium  faces.In Ex 4.3 we have
 H: a(x + y)   +z +t=0and slice has pentagon faces  In ex4.4 Webriefly  indicate special cases of  H: ax + by + z +t=0  with$ a > b >   $1.
As the methods in these examples is same asthe one in Ex 3.1  we   donot give  details.
In Ex.4.4 we consider the case of H:$ ax + by + cz +t=0$.Slice is  a paralleotope in case  $ a {\ge} b+ c + 1$  and  $b{\ge}  c + 1$.
If (i) and (ii) both fail then the slice has pentagon faces and is not a zonoid

We give these as  samples ; and do not consider every possible  case .Roughly  , the  non zonoids  prevail in our list of examples.

\textbf{ Diagram}They will help.

 Our \textbf{methods} are elementary   and can be found  for ex in [6]. We do use the   formula  for  vol of slices from [3] ( see also [10] ch1)  referred to above.
  
 We notethatthat  in al of our examples   we use the face [  t=  -1/2] of the cube , this is also a faceof the slice   The   $"domain"$  C of   face  is  found first , then an affine map T to determine theFace T(C). The points in C are found by  checking the x and y intercepts of lines involved   satisfy theconditions for slice:--$|x|$, $|y|$ and $|z|$ are all${\le}$ $1/2$. This condition must be satisfied by all coordinates of points  in the Face (of slices) that we  find,  and leads tothe conditions imposed on the coefficients of H.

Let us  first describe the result on slices  

.
\textbf{1.2Theorem  ([3][5] ,[7], [8] [13] ,[14] )}In all dimensions  the measure of  cube slice is between 1  and $\sqrt{2}$;  these arebest.

\textbf{ 1.3Zonoids}

Returning to the title of this  paper, about   Zonoids:--
Our  concern is :--“ When is a slice a zonoid?”   We  do  not have a  complete characterization of this .Instead let us concentrate  in  $R^4$, and give examples of  non zonoid  slices as well as those that are zonoids    and a consequence(  known) for $ I_4$  .
We  recall  from [ 6]   with X=$ R^ n$ :-- .  A \underline{ zonoid} is range of a  non atomic  vector measure and  above all  the classical Liapunov’s Theorem:- A  zonoid is compactand convex   A \underline{ zonotope} is sum of segments( each centered at the Origin) For our purpose we need the classic result of Herz and Lindenstrauss from [ 6]:-- The closed unit ball in every 2dimensional normed space   is a zonoid

\section{Zonoids and   Zonotopes}

\textbf{ 2.1 Theorem[ 6] }

i)If H is 2 dimensional then every such slice  is a zonotope
ii) In all dimensions every projection of  Q is a zonotope

\textbf{Proof:}

(i) This follows from the classic result  due to  Herz and Lindenstrauss quoted above and the result from ( [6] ) :-  in$ R^ 2$    every centrally symmetric polygon is  always  a   sum of segments

(ii) This is in [ 6] and can also be verified   directly. 
Hence the Theorem.

\textbf{Remark 2.2:} 

For much more about projections see, [ 4].

\section{ Example of slice that is not a  zonoid}

A  notedbefore,in contrast(Th 2.1, part i) to the situation in    $R^3$  we offer   an example of a slice in $R^4$  that is not a zonoid.  Reasons  to disqualify it from being a   zonoid  are  the  useful  facts, all from [ 6]  :--If K is a zonoid then

  (i) K has center of symmetry c say .In fact by definition of " K is  a Zonoid  "(see Introduction) 

  K= $\mu(\sum)$  for a ( vector measure) $\mu$  then  c= $ \dfrac{1}{2}\mu(S)$   will do   For, with $A^c$ = complement of set A, we have  
   $\dfrac{1}{2}(\mu(A) + \mu(A^c))$ =  $1/2 \mu(S) =c$ for every A in domain  $\sum $ of $\mu$

item{(ii) faces are translates of    zonoids of lower dimension and}
item  (iii) Since it has no  center of symmetry,the    triangle is  not a zonoid; neither is a trapezium   (trapezoid) or a pentagon

item(iv)Hence any compact ,convex, balanced set  that  has a triangular,( or a trapezium  face) cannot be a zonoid. Thus, the Octohedron   in $R^3$  is nota zonoid, for it   has triangular faces.  There are  deeper  non zonoids for ex the 1976 result due to LE Dor (  for ex[10]):--If $ 1<p <2$   and $ n{\ge3}$ then the closed unit Balls of the spaces  $l_n ^p$    are  not  zonoiids   

We give, in Th3.4,  a version of(ii) from [2 ]:-- a face   ( defined  below) is a translate     of some zonoid of lower dimension  . We need this version in the Th 3.4  to produce non zonoid  slices  in our examples. 

Let us recall  fom [6] theterm ,” \underline{Face} of   acompact convex set  K  “in a  real(  normed space ) X. Let us use “ H”  for  any hyperplane  ( not necessarily thru  O)

As above a \underline{hyperplane}   is
\begin{equation}\label{E:cong2}
   \ H =(x{\epsilon}X :( x ,x^* ) = {\alpha})\ ,
  \end{equation}
where $x^*$ isa non zero functional   in $ X^*$ and $ {\alpha} $is a real   number  .

The set  K is   “ on one side" of this H   if
\begin{equation}\label{E:cong3}
sup[( x, x*) : x \epsilon K]  \le\alpha,
\end{equation}

A similar     condition holds   with “ inf” replacing   
“ sup” and by$ \ge$ replacing $\le$ ;and\underline H supports  K  if  K is to one side  of H as in eq(3,)  H $\cap K  \neq \phi$  and K is not entirely  contained in H.  Finally the ( compact convex ) set H$ \cap$ K is called  \underline{“a Face  of K”} .  

Below we  use the fact    that an affine    map preserves  convexity. 

Let X,  Y be real Banach spaces  . Then a mapT:  X →  Y  is \underline{affine}  if    $ T ( ax  +  by)= aT(x) + bT(y) $   for every x, y in X and a , b $\ge  0 $ with a + b = 1.ie; the definition of  Linear map is now restricted to   line segments in domain.

Let us recall    K = $ H\cap Q $   is the slice corresponding to H[$ t=  -1/2$] In the next ( and other ) examples all we need is that the relevant  , y, z  coordinates of our points are limited  by  $|x|\le 1/2 $  etc . 

\textbf{3.1Example with triangle face}

Let us recall H  is given in $R^4$ by 
\begin{equation}\label{E:cong4}
x+ y+z + t=0  ,
\end{equation}
 
 The slice  (i) has triangular faces  and so is not a zonoid  

(ii) the intersections of slice with t= -c , $ 0 \le c <\dfrac{1}{2}$   
	
	are  hexagons ;these  are sections (iii) These   tend to the above  triangle as  c tends to  $1/2$
	
	\textbf { proof(i)} 
	
	Substituting   t  =  -1/2   in eq(4) of H , for any \textbf{x} =(  x, y, z, t)  in this H we have
	
	x  =  x( 1,0,0 -1)  + y (0,1,0,-1)  +z( 0,0, 1, -1)  is the  linear 
	combination x  u  +  y v  +  z w. (  these 3 vectors u, v, w are Linearly independent)
	
	First consider  the 2  dimensional  set   S in slice,in span of  vectors   u and v.  Starting with A  ( u/2)   on the “ x axis” and going counterclockwise,we see 
	
	that S is a   hexagon with vertices A ( u/2), B(v/2), C(( v-u)/2, A'= -A,

	B'= -B , C’= -C. Further  it is regular  all sides have  length $1/\sqrt{2}$ 
	and that this = sum of 3 segments , OA  , OC  and OB'. This  set   S is in plane   z=0
		Now we consider the  3rd term in above eq   for x  ; we   note  that   the vector  D=  w/2   cannot be added to   A or B   as the sum will leave the cube
		We consider  the Hyperplane 
	 $H_1$  = { ( x, y, z, t):   t= -1/2} )   or,  simply  by t  =  -1/2  and claim that 
		
		( a)  this  plane supports   the slice K   and that 
		
		( b)the face F  =$ H_1{\cap} K $  is  convex triangle.
		
		As  noted above, ( b) disqualifies the slice from being   a zonoid

		Let us verify the  claims. Now (a) follows      directly from def. of Q.In fact	for every  element in Q we have $t{\ge}-1/2 $ ie.,  Q  is to “ one side “ of $H_1$; 	so is the slice .Further , the elements  A , B,    are in the slice, and   also  lie in$ H _1$ , hence  in Face F.  The Origin O  is in slice K  not in $H_1$; ie. the sliceis not  entirely  contained  in$ H_1 $.Hence  $H_1$   is a  supporting hyperplane of the slice as claimed.

		For claim (b)   we may write any element    in the Face  as 
		
		\textbf{f (x) = (x, y, z, t)=( x, y, 1/2 -  x- y , - 1/2})    , 
	
	  since  we use  t=  - 1/2  in eq   (4) of H and
	
	we  get z=  $ 1/2--x-y $.
	
	As x is in Q we need $| x|$ and  $| y|$ and also $ |z|$   from above  ${\le}1/2$ and so  
	
	\begin{equation}\label{E:cong5}
	   |1/2- x- y| \le 1/2 
	\end{equation}
	
	This last translates to
	
	\begin{equation}\label{E:cong6}
         0 \le  x+y  \le  1  
	\end{equation}
	
	Geometrically,   we note that the   last inequality   gives  two boundary lines of $"domain C"$ say$L_1$ := x+y =1 , and  $L_2$:=  x+ y=0.
	We sketch these lines; as the x-intercept of$L_1$  exceed the bound 1/2  let us consider its intersection with the line x =1/2 toget point  (1/2, 1/2).; intersection of  $L_2$ with the line  y=1/2   gives ( -1/2, 1/2). This line with  y=-1/2  gives (1/2, -1/2)
	
	 These  result in a ( convex right angled ) triangle C in x-y plane  with   above
	
	vertices     P( 1/2, -1/2) , Q( 1/2,  1/2)  and R ( - 1/2, 1/2)

	Now let us define a map Tfrom Cto F by 

	\begin{equation}\label{E:cong7}
	\ T ( x,y) = ( x, y,    1/2   -x-y,   - 1/2)\,
	 \end{equation}
	 
    and  C is its domain. Then   we may verify that,  T is affine and that T( C)=F.
	
	Further   as   observed before statement  of this example,  affine  map
	
	preserves convexity,   and so the image T ( C) =  convex hull of the 3 points
	
	($ p_1$,$ p_2$   , $p_3$) where $p_1$   = T( P)=  ( 1 /2  ,  -1/2 ,1/2,  -1/2) , $ p_2 $ = 
	
	T( Q)= ( 1 /2,  1/2, -1/2, -1/2) and$ p_3 $= T( R)= (- 1 /2,  ,1/2,  1 /2,-1/2).   These  points   are not collinear, form a triangle and we conclude that the face F  is a triangle , completing Claim ( b) and proof of (i) .
	
		We need to prove   (ii)  and (iii) .   

     Recall K = slice  ; now  we    let	$ 0 < c <  1 / 2$  and   Section  $K_c = K \cap[t  =  -  c]$.   
	
	Use  t=  - c in   eq  (4) H;   any x   in   Kc  is   then   of 
	
	the form      \textbf{x} =   ( x, y , -x-y+c,    -c)  with the conditions
	
	  $|x|$, $|y|$ and $|x + y -c| \le  1 /2$   .
	
	Similarly to above ( 6) this last translates to 
	
	\begin{equation}\label{E:cong8}
         \ -1/2 +c {\le} x+y {\le}1 / 2  + c \ ,
	  \end{equation}
	
		As in part (i)  we   draw the $" boundary"$  lines  $L_1$  , $L_2$  from eq (8)  .  Again , both  the  x and y- intercepts of $L_1$  fail the  bounds   of 1/2; however $L_2$  passes ( noting the limits on c)Then  we  find the  vertices  of our    $" domain C"$  by intersecting   $L_1$  and $L_2$   with the  lines y=1/2, y=  -1/2,   x=  1/2  and x=  - 1/2. . We  get a hexagon( domain). Its  6  vertices are shown in a Chart in next 	Theorem  3.3  and as follows:--
	
	$p_6$  =(c, -1/2)  on lines  y =  -1/2   and $L_2$  , $ p_1$  = 
	
	( 1 / 2, -1/2)  and$ p_2 $ =( $1/2$  , c)  on line  x=   1 / 2 and$L_1$  .Next $p_3$  =( c, 1 / 2) on lines $L_1$  and y=1/2 and  $p_4 $ =(- 1/2,  1/2)  then $p_5$  =  ( -1/2, c)  on    lines $L_2$  and  x=  -1/2

	These   6   points $( p_i  )$ form a   hexagon  making    
	
	the new domain   “C”  of  map T defined  analogous to eq(7) in part (i) above.
	
	As there we see that the  section  $K_c$  = T( C ) is also a hexagon.
	
	Finally  let  c tend to 1  / 2  ; then we  see from   above that the  
	
	following    vertices   coincide:-  $ p_2 $ = $p_3$  =(1/2, 1/2)  ,$ p_4$ = $ p_5 $
	
	=(-1/2,  1  /2)  and   $ p_6$ = (1/2, -1/2) = $ p_1$  . Correspondingly ( as in 
	
	case i above)  we verify that  the   section T( C ) becomes the   triangle  in 
	
	part(i)  completing thereby proof   of (ii) and  the example

	\textbf{ Remark 3.2.}   Above we  used    the  hyper  plane  given by the equation,   
	
	t=  -   1  /2 and       found that  the face of slice    given  by it  is  triangular 
	
	; we  may instead consider  t = 1  / 2  Further, the equation defining H is 
	
	symmetric with respect to the four variables   x, y, z, t. Hence  we  may conclude 	that there are 8  triangular faces.  We do  not     know   what are the remaining 	faces   and we  think there are  4 more  but not triangles  .

	For the next result,   we  follow [  3]Lemma 3  ( see also[10] ch1  ) and  recall from Introduction eq(1) theintegral  $I_p$ :
	
$\dfrac{1}{\pi}\int_{R}\dfrac{|sin t| ^p}{|t| ^p}  dt $
.

Here  p  is an integer $\ge$2, and we have from the result in [3] above,

  the formula forthe  exact value of slice :-     

\begin{equation}\label{E:cong9}
      \|H{\cap}Q|  = \dfrac{1}{\pi} \int_{R} g(t) dt\\ ,
 \end{equation}     

 where g is the   finite product

\begin{equation}\label{E:cong10}
     \  g(t) = \prod_{1}^{N} \dfrac{sin a_it}{a_it}\\ ,
\end{equation}

	and the sequence    ( $ a_i $ ) ( of coordinates  of vector normal to H) is  normalized in $l ^2$  and   also each $ | a_j |  \le 1/ \sqrt {2} $
	
	To  find volume of   slice  S  we use Cavaleri  ‘s principle = Fubini's Theorem . 
	
	Let    $| A(  c) |$=  area    of the section  of  S  by  plane [  t= -c]. Then  the 
	
	vol of  slice  =2  $\int _{0}^{ 1/2} |A( c) | dc$  . 
	
	We saw in   ex 3.1  that A (  c)  is a hexagon .  We give the  details  in the 
	
	next result;

	\textbf{ 3.3Theorem} (i)The volume of the slice  in Ex3.1  is   4/3   (ii) $I_4 $ =   2/3

	Proof(i)    We refer to  part (b)in ex3.1 and list the    vertices of the  hexagons in 
	
	domain C  as well as in  the range  T(C).

	 Recall  T( x, y)= ( x, y, c-x-y, -c)   	with $ 0  < c <  1  /2 .$

	Domain C ............................... rangeT( C)
	
	$ p_1(1/2, -1/2)$...................      $P_1(1/2, -1/2, c, -c)$

$ p_2$ (1/2,  c) ........................   $P_2 $ (1/2,  c ,-1/2 , -c)

$ p_3$ ( c , 1/2)  ....................       $P_3 ( c , 1/2,   -1/2, -c)$

$p_4$  ( -1/2, 1/2) ................   .. $P_4(-1/2, 1/2, c, -c)$

$p_ (-1/2,  c)$  ......................
 $P_5(-1/2, c, 1/2, -c)$

$p_6( c, -1/2)$ ..........................$P_6 ( c, -1/2, 1/2 , -c)$

We  claim that  area of  domain C=

\begin{equation}\label{E:cong11}
   | A( c)| =  (3 /4  - c ^2 ),
\end{equation}

In the following we use  formula     for area of trapezium  by  rule

(1 /2) h (a + b) where h is the height  and a, b    are  lengths of parallel sides. 

Let us use the  chart for domain C  first then use it to get the   image.

The  domain C = two trapeziums$ T_1$ and$ T_2 $ ;  these are the top and  at 

bottom  respy.  Namely, $T_1$  has vertices, $p_5$, $p_2$, $p_3$ and $p_4$ 

 and $T_2$  has vertices, $p_6$, $p_1$,$p_2$ and$p_5$.

 Then we have

$| T_1 |$= $ \dfrac{1}{2}( 1+c +1/2) ( 1/2  - c)= 1 /2 ( 3/2 +c) (  1 /2  -c)$    and

$|T_2 |$=  $ 1/2 ( 1  -c +1/2) ( 1/2  + c)=   1 /2(  3/2   -c) ( 1/2   +c) $ 

Adding them we get the    eq ( 11) for  A(c).

To get the  area of  image T( C )   observe that the domain C  is the 

projection on plane(  z=0) of the  wanted T( C) .

For the factor   needed we   note  that the unit        normal to H is

n=( 1/2 , 1/2 , 1 /2 ,  1 /2) and that $ e_3 $  =    ( 0,0,1, 0).

Using the dot product $ n . e_3 $ we see that area of     Projection

= 1/2  area of T (  C). Thus the area of T (C) = $2( 3 / 4-  c^2 )$  from above. 

We  integrate from c=0 to 1 / 2   to get  

 $2\int_{0}^{1/2}(\dfrac{ 3 }{  4}   -   c  ^2  ) dc=  2/3 $

Taking into account alsothe part t=  1/2  to 0 we      get  2( 2/3)  =4/3 as claimed

Part (ii)  :Recalling H:  x+ y + z+t =0  and  the  coefficients, normalized , we   apply the formula from [  3] quotedabove   in  eq   (9) , (10)  to get   vol of slice=

$\dfrac{1}{\pi} \int_{R}  \frac { (sint/2) ^4}{(t/2)^4}  dt  $  .

(  as in part (i)    we used each $a_i $ coefficient  is 1/2     due to normalizing them  in eq of H )    Now   a change of variable gives  

$\dfrac{2}{\pi}\int_{R}$ $ ( sint/t ) ^4 $   dt  = $2 I_4 $.

From part (i) we have  $ 2I_4$ =  4/ 3 

and so part (ii)  and the Theorem

Above,in example  of a  non zonoid slice we used  the important  fact about faces of a zonoid from [6] (  innext Theorem)  The  following proof is different  from the one  in[6] which uses “ Every Zonoid  is a zonoid of  moments” .This  approach  is  not suitable for  our  purpose; hence we give a proof ( in[ 2]) in next result.  We see in  the proof  that it is more meaningful incase the Face isnot a singleton, ie.   when the composed    measure $x ^*o{\mu}$  is not equivalent to  ${\mu}$

\textbf{3.4Theorem[6][2]]}Let K =$ {\mu}({\sum} ) $ be a    zonoid in  X=$ R^n$, and  H a  supporting Hyperplane given by x* in X*.Then   the face F=$ K {\cap}  H $ is a translate of a zonoid of   lower dimension . In fact there  are $  {\mu}$ almost  disjoint  sets  $ S_0$  and $S_1$   such that

(i)$ x  ^*o{\mu}(S_1)$=sup$\{ x ^*o{\mu}( E): E {\epsilon} {\sum }\}$   and every set  E in$ S_0$ is  $ x ^*o {\mu} $ - null

(ii F= $\mu(S_ {1})   +  \mu _{S _0} (\sum )$.

\textbf{ Proof}: With $\beta  =  supx^*\mu(\sum)$ we have,from definition of   F 

\begin{equation}\label{E:cong12}
      \  F= \{ x:x={\mu}( E)    s. t.     x^*(x)={\beta}\} \\  ,
 \end{equation}

Let $S^+$  be such that   $x^*o\mu(S^+)=\beta  $. 

We will write  $ S^+  = S_0  \cup  S_ 1$  as stated in    the Theorem.

To  do this let usnote  that the  signed measure  $x^*o\mu  <<\mu $

;consider those  E  that are   $ x^*o \mu$ - null but not $\mu$  - null.   

(ifthere are no such  sets E then $ S_0$may be taken to be $\emptyset$ ).

Othewise consider a maximal pairwise disjoint family    ofsuch sets; this  family is countable, so that theirunion is in   $\sum$.
Call this set $S_0$  and let $S_1$  = $S ^+  -  S_0$  .Then 

(i) follows from the fact that 

$x ^* o \mu (S_0 )$ =0 by the    construction of $S_0$     and so

$\beta  = x ^*  o \mu(S^+ )  $=  
 $ x  ^* o\mu (S_0  )  + x ^*o \mu  ( S_1 )$

= $ x ^* $o  ${\mu }  ( S_1  )$ . we see that     second part in (i) follows by  construction again.

As for  part(ii) we  have from Eq (12)  if   x= ${\mu}( E)$    ${\epsilon}$ F then 

$x^*o {\mu}$ ( E) = ${\beta}$ .

We need to write x =$ {\mu}$ ( E) asthe sum, ${\mu}$ ( E)= $ {\mu}$  ($S_1 $)  +  ${\mu}$ ( A)  for some set  A$\subset  S_0$
To do this,first we claim that( ae  $|x^*o\mu  |$  ) this  $E{\subset} S^+$  . If not  we can argue to   contradict to  the fact that S  = $S^+ \cup S^- $  is a  Hahn  decomposition of the underlying set S  in terms of  $x^*o\mu$      

Again we can argue   that  $S_1$ - E is ${\mu}$– null;from it being $x^*o{\mu}$ null, and then on ( subsets of ) $S_1$     these two  measures  are  equivalent by construction.

Hence we have E= $ E \cap   S_ 1 \cup  E \cap S_0 $, and so

$\mu(  E) =  \mu ( E \cap) S_ 1 )  +  \mu   (E \cap )S_0 )$ =

$\mu(S_ 1 )  + \mu  (A) $ with A = $ ( E \cap S_0) $         $ \subset  S_0 $ as claimed.

Hence the Theorem

\section { Examples of non zonoids with  pentagon faces  and some zonoids  a $>$1 }.

As in Introduction welet H: ax + by+ cz + t=0 be a  hyperplane in $R^4$,with 
 $a \ge  b \ge c \ge  1  $.

 We donot  consider all cases   but hope the following are of interest. There are non trivial cases  of zonoid slices .   As the  methods  are same as the  one in the earlier ex 3.1 we only    summarise the results  It seems the non zonoid  slices  dominate:--
 
 In the next  ex. we donot know if the converse is true in this generality. Hencewe  give some special cases of the eq of H in the EXs 4.2 and on.In all cases  for  the Face    we use as before  the support hyperplane of 
 Q  $ [ t= \dfrac{-1}{2}$]
 
 \textbf{4.1 H: general caseabove } 
 
 If   $ a  \ge  b+c+1$  then  the sliceis a zonotope.

Proceeding as in Ex3.1, we find  the " domain" for the face .For this we have the boundary lines $L_1$ to be ax + by = $\dfrac{c +1}{2}$ and $L_2$ to be ax+by =$\dfrac{-c +1}{2}$

  Firstwe note    bothx andy-intercepts of $L_2$    are always ( regardless ofthis condition ) $\le \dfrac{1}{2}$ in absolute value. As  for$L_1$  this condition  gives the x-intercepttobe ${\le}  \dfrac{1}{2}$  in absolute velue.   In the   following " domain  the vertex $p_2$ depends on this condition, ie. its " $|x|$  satisfies the limits $\le \dfrac{1}{2}$. 

 Withthe condition above we have  now the chart

 \textbf{Domain}           
 
 $p_1$ ($\dfrac{b-c+1}{2a}$ , $\dfrac{-1}{2}$ )
 		
 $p_2$  ($\dfrac{c+1+b}{2a}$ , $\dfrac{-1}{2}$  )	
 
 $p_3$  ($\frac{c+1-b}{2a}$,  $\dfrac{1}{2}$  ) 
 
 $p_4$   ($\dfrac{1-b-c}{2a} $ , $\dfrac{1}{2} $ )

 Next the corresponding points on the Face:--

 \textbf{Face}T(x,y)= $(x, y,\dfrac{1/2-ax-by}{c}, \dfrac{-1}{2}  )$

 $P_1$ ($\dfrac{b-c+1}{2a}$, $\dfrac{-1}{2}$, $\dfrac{1}{2}$, $\dfrac{-1}{2} $)
 
 $P_2$  ($\dfrac{c+1+b}{2a}$,  $\dfrac{-1}{ 2}$, $\dfrac{-1}{2}$, $\dfrac{-1}{2}$)
 
 $P_3$   ($\dfrac{c+1-b}{2a}$ , $\dfrac{1}{2}$, $\dfrac{-1}{2}$ ,$\dfrac{-1}{2}$)
 
 $P_4$  ( $\dfrac{1-b-c}{2a}$,  $\dfrac{1}{2} $, $\dfrac{1}{2}$,$\dfrac{-1}{2}$ )

  It is seen   that this domain is a parallogram with the parallel sides( so is the Face): 
 
 ($p_1$$p_2$ )  = ($p_4$$p_3$) =  ($ \dfrac{c}{a} $,  0)  and
 
 ($p_2 p_3$)  = $ (p_1$ $p_4$)   =($\dfrac{b}{a} . -1)$
 
 Likewise,  it   can be verified  using the map "T"  , that so is  the  Face.

 Further, the sections of thesliceby planes with eqs   t= $-c_1$, with $ 0 < \textbf{}c_1  < 1/2 $   are  parallograms that are congruent to the one for the   Face. Hence  it follows (using symmetry)  thatthe slice is a zonotope.

\textbf{  4.2 H : a x +y + z+ t=0} 

In one  direction  this is a special caseofEx 4.1  ;however  due to limitation of eq ofH   we can state " iff"   and we  give  details :--

In this case if $ a \ge3 $ then the  slice isa zonoid ; it is a   paralleotope if not  the slice has    pentagon faces and is   not a zonoid.

\underline{ Case  $ a \ge3$:}  Face  is a ;paralleogram ; so is   every parallel section  congruent to it

Letus note that analogously to  eq(6) above we replace  x  by ax  there. Thus the x-intercept of the  line  with equation $ ax +y =1$  is   x= $1/a$.The condition $ x {\le}$  $ 1/2$   now holds (due tothe condition on a ). This forces the "Domain" to be a paralleogram as we now state. As  before  we use  (x,y) forpoints $p_i$     and 

 \textbf{T( x,y)= (x, y, 1/2- ( ax +y),  -1/2)} for points $P_i$  : 
      Domain  C(x,y)             Face  T(C)
        
    $ p_1$( 1/2a,  -1/2)                         $P_1$(1/2a, - 1/2, 1/2, -1/2)

     $p_2$(3/2a, -1/2)                            $P_2$ ( 3/2a,  -1/2,  -1/2, -1/2)

     $p_3$ (1/2a, 1/2)                           $P_3$(1/2a, 1/2, - 1/2,  -1/2)

     $p_4$ (-1/2a, 1/2)                           $ P_4$(- 1/2a, 1/2, 1/2, -1/2) 
       
     We see that the  opposite sides are  parallel and  have equal     length , so that the Face is a rhombus .Further so  is any section   by plane[ t= -c] with $ 0< c < 1/2$, the  area does not depend on     c and  equals $\sqrt{1 + 2 a^-2}$

     \underline{ case  $a  < 3$}  in this case we can verify  the " domain" tobe apentagon; so is the face and  slice is not a zoniod

    \textbf{4.3  H: a(x +y) +z +t  =0}  with $ a{\ge} 2$
    
    The Face$ [ t=- 1/2]$  is a trapezium  again,  slice not a zonoid
    
     \textbf{ 4.4 H:ax +by+ z +t =0  }( compare  ex4.1 ) 
     Face is a paralleogram  in case $ b+2{\le}a  $  . The parallel sections$[t= -c]$ are  congruent parallograms, and the slice isa parallotope.   Otherwise Face is a pentagon,    slice is not a zonoid
     
     \textbf{4.5  H: ax + by + z +t =0}
     The slice is a zonotope if (i)  a ${\ge} b + c+ 1$ and (ii) $b{\ge} c+ 1$.
     
     In case (i)  and (ii)  both  fail   Face is a pentagon and slice is nota zonoid.
     
     If  (i) fails  but(ii) is true then the Face is a hexagon
     
     \textbf{Remark 4.5} In the   last case we dont know if the slice is a zonoid

\textbf{ADDED   REFERENCE JULY 2023}   Ethan Bolkerin  2021  Joint Mathematics Meetings(//exhibitions//2021-jointmathematicsmeetings-meetings)

\textbf{ Acknowledgements}: 

The author is  thankful to an anonymous referee( of another      related paper)  for the reference [ 12]  and the  information there about the Sinc Integral $I_4$

Further  ,he is most grateful to Professor Ethan Bolker  for  the reference above  in Joint Math Mtgs 2021

\textbf{Announcement}The  abstract of this paper wasannounced in the Proceedings of Ramanujan Birthday Conference   held on Dec 22nd 2017 in IITMadras, CHennai, India

	date first done : Dec 20th 2016
	
	Retired from
	SUNY/College @Old     Westbury
	Old Westbury NY  11568-0210
	e mail  rajan.anantharaman@gmail.com
	
	address: 3032 Briar Oak Dr 
	         Duluth GA 30096

	MR   Classification:
	Primary  52 A 20,  Secondary  42 A38,  52A40
\end{document}